\documentstyle[12pt,amssymb,amscd]{amsart}  

\newfont{\cyr}{wncyr10}

\newtheorem{thm}{Theorem}

\newtheorem{lemma}[thm]{Lemma}

\newtheorem{proposition}[thm]{Proposition}

\theoremstyle{remark}
\newtheorem{remark}[thm]{Remark}

\theoremstyle{definition}
\newtheorem{definition}[thm]{Definition}

\numberwithin{equation}{subsection}

\newcommand{\Ai}{{A_{\infty}}}
\newcommand{\HCp}{\operatorname{HC}_{\bullet}^{\operatorname{per}}}
\newcommand{\HCn}{{\operatorname{HC}}_{\bullet}^-}
\newcommand{\CCp}{\operatorname{CC}_{\bullet}^{\operatorname{per}}}
\newcommand{\CCn}{{\operatorname{CC}}_{\bullet}^-}

\newcommand{\ChA}{C_{\bullet}(A)}
\newcommand{\CochA}{C^{\bullet}(A)}

\newcommand{\CCnCocha}{{\operatorname{CC}}_{\bullet}^-(C^{\bullet}(A))}
\newcommand{\HCnA}{\HCn(A)}
\newcommand{\CCnA}{\CCn(A)}

\newcommand{\CCpA}{\CCp(A)}

\newcommand{\gA}{{\mathfrak{g}}^{\bullet}_A}

\newcommand{\Ccc}{C^{\bullet}(A)}

\newcommand{\Hcc}{H^{\bullet}(A)}

\newcommand{\CCc}{C_{\bullet}(C^{\bullet}(A))}

\begin{document}

\title{On the Gauss-Manin connection in cyclic homology}

\author{Boris Tsygan}
\email{tsygan@@math.northwestern.edu}
\address{Northwestern University}

\maketitle

\smallskip

{\centerline{In memory of Yu.~L.~Daletsky}}

\smallskip

\section{Introduction} For an algebraic variety $S$ over a commutative field $k$ of characteristic zero, let $A$ be a locally free ${\mathcal O}_S$-module which is an associative ${\mathcal O}_S$-algebra. In \cite{G}, Getzler constructed a flat connection in the ${\mathcal O}_S$-module $\HCp(A)$, the periodic cyclic homology of $A$ over the ring of scalars ${\mathcal O}_S$. This connection is called the Gauss-Manin connection. In this paper we define this connection at the level of the periodic cyclic chain complex $\CCp(A)$.

Recall that for an associative algebra over a commutative unital ring $K$ one can define the Hochschild chain complex $\ChA$, the negative cyclic complex $\CCn (A)$, and the periodic cyclic complex $\CCp (A)$, as well as the Hochschild cochain complex $\CochA$ (\cite{L}, \cite{T}, \cite{Ge}). The latter is a differential graded Lie algebra, or a DGLA, if one shifts the degree by one: ${\frak{g}}^{\bullet}_A=C^{\bullet +1}(A)$. Recall that $\CCnA =(\ChA[[u]], b+uB)$ is a complex of $K[[u]]$-modules. Here $u$ is a formal variable of degree $-2$. We can view $\CCnA$ as a {\em cochain} complex if we reverse the grading. In particular, the {\em cohomological} degree of $u$ is $2$. The complex $\CCnA$ is known to be a DG module over the DGLA $\gA$, the action of a cochain $D$ given by the standard operator $L_D$ (cf. \cite{T} or \ref{eq: L} below).

Consider another formal variable, $\epsilon$, of degree $1$. Now consider the DGLA 
\begin{equation} \label{eq:gA(u,e)}
(\gA[u, \epsilon], \delta+u\frac{\partial}{\partial \epsilon})
\end{equation}
\begin{thm} \label{thm:linfty structure}
On $\CCnA$, there is a natural structure of an $L_{\infty}$ module over $(\gA[u, \epsilon], \delta+u\frac{\partial}{\partial \epsilon})$. This structure is $K[[u]]$-linear and $(u)$-adically continuous. The induced structure of an $L_{\infty}$ module over $\gA$ is the standard one.
\end{thm}
We recall that an $L_{\infty}$ module structure, or, which is the same, an $L_{\infty}$ morphism $\gA[u,\epsilon] \to \operatorname{End}_{K[[u]]}(\CCnA)$, can be defined in two equivalent ways. One definition expresses it as a sequence of DGLA morphisms 
\begin{equation}\label{eq:L-infty as domik}
\gA[u,\epsilon]\leftarrow {\mathcal L}\to \operatorname{End}_{K[[u]]}(\CCnA)
\end{equation}
where the morphism on the left is a quasi-isomorphism. Alternatively, one can define this $L_{\infty}$ morphism
as a collection of $K[[u]]$-linear maps 
\begin{equation} \label{eq:linfty-morphism}
\phi _n: S^n( \gA[u,\epsilon] [1])\to \operatorname{End}_{K[[u]]}(\CCnA)[1]
\end{equation}
satisfying certain quadratic equations. Using that, one can define Getzler's Gauss-Manin connection at the level of chains as a morphism
$$
\Omega ^{\bullet}(S, \CCpA)\to \Omega ^{\bullet}(S, \CCpA)
$$
of total degree one such that
$$\omega \mapsto d\omega +\sum _{n=1}^{\infty} \frac{u^{-n}}{n!}\phi _n(\theta, \ldots, \theta)
$$
where $\theta$ is the $\gA[u,\epsilon]$-valued one-form on $S$ given by 
$$\theta (X)(s) = L_X m_s \epsilon.$$
Here $s\in S$, $m_s$ is the multiplication on the fiber $A_s$ of $A$ at the point $s$, and $X$ is a tangent vector to $S$ at $s$.

A few words about the proof of the main theorem. We define the $L_{\infty}$ morphism by explicit formulas (Theorem \ref{thm:action of twisted bar} and Lemma \ref{quis of twisted and untwisted}), but the proof that they do satisfy $L_{\infty}$ axioms is somewhat roundabout. Recall that the Hochschild cochain complex $\CochA$, with the cup product, is a differential graded algebra (DGA). One can consider the negative cyclic complex $\CCnCocha$ of this DGA. In \cite{TT} and \cite{T}, an $A_{\infty}$ structure on this complex is constructed. The negative cyclic complex $\CCnA$ is an $A_{\infty}$ module over this $A_{\infty}$ algebra. From this, we deduce that $\CCnA$ is a DG module over some DGA which is related to the universal enveloping algebra $U(g_A [u, \epsilon])$ by a simple chain of quasi-isomorphisms.

A statement close to Theorem \ref{thm:linfty structure} was proven in \cite{DT}. Our proof substantially simplifies the proof given there. Note that a much stronger statement can be proven. Namely, $\CochA$ is in fact a $G_{\infty}$ algebra in the sense of Getzler-Jones \cite{GJ1} whose underlying $L_{\infty}$ algebra is $\gA$, (\cite{Ta}, \cite{H}); moreover, the pair ($\CochA$, $\ChA$) is a {\it homotopy calculus}, or a $\operatorname{Calc}_{\infty}$ algebra (\cite{KS}, \cite{TT}, \cite{T}). The underlying $L_{\infty}$ module structure on $\CCnA$ is the standard one. From this,  Theorem \ref{thm:linfty structure} follows immediately. (The interpretation of the $A_{\infty}$ algebra $\CCnCocha$ in terms of the $\operatorname{Calc}_{\infty}$ structure is given in \cite{TT}). However, theorems from \cite{KS}, \cite{TT}, \cite{T} are extremely inexplicit and the constructions are not canonical, i.e. dependent on a choice of a Drinfeld associator. Our construction here is much more canonical and explicit, though still not perfect in that regard. It does not depend on an associator; it provides an explicit structure of a DG module over an auxiliary DG algebra (denoted in this paper by $B^{\operatorname{tw}}(\gA[u, \epsilon])$). Unfortunately, this auxiliary DGA is regated to our DGLA somewhat inexplicitly.

Theorem \ref{thm:linfty structure} implies the existence on $\CCnA$ of a structure of an $A_{\infty}$ module over $U(\gA [u, \epsilon])$; the induced $A_{\infty}$ module structure over $U(\gA)$ is defined by the standard operators $L_D$. An explicit linear map 
$$U(\gA [u, \epsilon])\otimes _{U(\gA)}\CCnA \to \CCnA$$
was defined in \cite{NT}. It is likely to coincide with the first term of the above  $A_{\infty}$ module structure. 

This paper is, to a large extent, an effort to clarify and streamline our work \cite{DT} with Yu.~L.~Daletsky. I greatly benefited from conversations with P.~Bressler, K.~Costello, V.~Dolgushev, E.~Getzler, M.~Kontsevich, Y.~Soibelman, and D.~Tamarkin.

This work was partially supported by an NSF grant.

\subsection{The Hochschild cochain complex} \label{hocochain}

Let $A$ be a graded algebra with unit over a commutative unital ring
$K$ of characteristic zero.  A Hochschild $d$-cochain is a linear map $A^{\otimes d}\to A$.  Put,
for $d\geq 0$,
$$
        C^d(A) = C^d (A,A) =\operatorname{Hom}_K({\overline{A}}^{\otimes d},A)
$$
where ${\overline{A}}=A/K\cdot 1$. Put
$$
        |D|\;=({\em {degree\; of\; the\; linear\; map\; }}D)+d  
$$
  Put for cochains $D$ and $E$ from $C^{\bullet}(A,A)$
$$
        (D\smile E)(a_1,\dots,a_{d+e})=(-1)^{| E|\sum_{i \leq d}(|a_i| + 1)}
        D(a_1,\dots,a_d) E(a_{d+1},\dots,a_{d+e});
$$
$$
        (D\circ E)(a_1,\dots,a_{d+e-1})=\sum_{j \geq 0}
        (-1)^{(|E|+1)\sum_{i=1}^{j}(|a_i|+1)}
D(a_1,\dots,a_j,  
        E(a_{j+1},\dots,a_{j+e}),\dots); 
$$
$$
        [D, \; E]= D\circ E - (-1)^{(|D|+1)(|E|+1)}E\circ D
$$
These operations define the graded associative algebra
$(C^{\bullet}(A,A)\;,\smile)$ and the graded Lie algebra
($C^{\bullet + 1}(A,A)$, $[\;,\;]$) (cf. \cite{CE}; \cite{G}).
Let
$$
        m(a_1,a_2)=(-1)^{\deg a_1}\;a_1 a_2;
$$
this is a 2-cochain of $A$ (not in $C^2$).  Put
$$
        \delta D=[m,D];  
$$
$$
        (\delta D)(a_1,\dots,a_{d+1})=(-1)^{|a_1||D|+|D|+1} 
 a_1 D(a_2,\dots,a_{d+1})+  
$$
$$
        +\sum
_{j=1}^{d}(-1)^{|D|+1+\sum_{i=1}^{j}(|a_i|+1)}
                D(a_1,\dots,a_ja_{j+1},\dots,a_{d+1})  
 +(-1)^{|D|\sum_{i=1}^{d }(|a_i|+1)}D(a_1,\dots,a_d)a_{d+1}
$$
 
One has
$$
        \delta^2=0;\quad\delta(D\smile E)=\delta D\smile E+(-1)^{|deg D|}
                D\smile\delta E  
$$
$$
        \delta[D,E]=[\delta D,E]+(-1)^{|D|+1}\;[D,\delta E]
$$
($\delta^2=0$ follows from $[m,m]=0$).

Thus $ C^{\bullet}(A,A)$ becomes a complex; we will denote it also by $C^{\bullet}(A)$. The cohomology of this complex 
is $H^{\bullet}(A,A)$ or the Hochschild cohomology. We denote it also by 
$H^{\bullet}(A) $.  The $\smile$ product induces the
Yoneda product on $H^{\bullet}(A,A)=Ext_{A\otimes A^0}^{\bullet}(A,A)$.  The operation
$[\;,\;]$ is the Gerstenhaber bracket \cite{Ge}. 

If $(A, \;\; \partial)$ is a differential graded algebra then one can define 
the differential $\partial$ acting on $A$ by 
$$
\partial D \;\; = \; [\partial , D]
$$

\begin{thm} \cite{Ge} The cup product and the Gerstenhaber bracket induce a Gerstenhaber algebra structure on $\Hcc$.
\end{thm}
For  cochains $D$ and $D_i$ define a new Hochschild cochain by the following formula of Gerstenhaber (\cite{Ge}) and Getzler (\cite{G}):

$$D_0\{D_1, \ldots , D_m\}(a_1, \ldots, a_n) = 
\sum (-1)^{ \sum_{k\leq {i_p}}(|a_k| + 1)(| D_p|+1)}  D_0(a_1, \ldots ,$$
$$a_{i_1} , D_1 (a_{i_1 + 1}, \ldots ),\ldots ,
D_m (a_{i_m + 1}, \ldots ) , \ldots)$$
\begin{proposition}\label{prop:brace structure}
 One has
$$
(D\{E_1, \ldots , E_k \})\{F_1, \ldots, F_l \}=\sum (-1)^{\sum _{q \leq i_p}(|E_p|+1)(|F_q|+1)} \times 
$$
$$
\times D\{F_1, \ldots , E_1 \{F_{i_1 +1}, \ldots , \} , \ldots ,  E_k \{F_{i_k +1}, \ldots , \}, \ldots, \}
$$
\end{proposition}

The above proposition can be restated as follows.  
For a cochain $D$ let $D^{(k)}$ be the following $k$-cochain of the DGA $\CochA$:
$$
D^{(k)}(D_1, \ldots, D_k) = D\{D_1, \ldots, D_k\}
$$
\begin{proposition} \label{etoee}
 The map 
$$
D \mapsto \sum_{k \geq 0} D^{(k)}
$$
is a morphism of differential graded algebras
$$ C^{\bullet}(A) \rightarrow C^{\bullet}(  C^{\bullet}(A))$$
\end{proposition}
\subsection{Hochschild chains}  \label{ss:hochschild-1}
Let $A$ be an associative unital dg algebra over a ground ring $K$. The differential on $A$ is denoted by $\delta$. Recall that by definition
$$\overline{A} = A / K\cdot 1$$
Set
$$C_p (A,A) = C_p(A) = A \otimes \overline{A} ^{\otimes p}$$
Define the differentials $\delta: C_{\bullet}(A) \to C_{\bullet}(A)$, $b: C_{\bullet}(A) \to C_{\bullet - 1}(A)$, $B: C_{\bullet}(A) \to C_{\bullet + 1}(A)$ as follows.
\[
\delta (a_0\otimes\cdots\otimes a_p ) = 
\sum_{i=1}^p {(-1)^{\sum_{k<i}{(| a_k| + 1)+1}}
(a_0\otimes\cdots\otimes\delta a_i \otimes \cdots \otimes a_p )};
\]

\begin{equation} \label{eq:b grad}
b(a_0 \otimes \ldots \otimes a_p) = \sum _{k=0}^{p-1} (-1)^{\sum_{i=0}^{k} {(|a_i| + 1)+1}}
 a_0 \ldots \otimes a_k a_{k+1} \otimes \ldots a_p 
\end{equation}
$$+ (-1)^{|a_p| + (|a_p|+1)\sum_{i=0}^{p-1}(|a_i|+1)} a_pa_0 \otimes \ldots \otimes a_{p-1};
$$

\begin{equation} \label{eq: B graded}
B(a_0 \otimes \ldots \otimes a_p) = \sum_{k=0}^p (-1)^{\sum _{i \leq k}(|a_i| + 1) \sum _{i \geq k}(|a_i| + 1)} 1 \otimes a_{k+1} \otimes \ldots a_p \otimes
a_0 \otimes \ldots \otimes a_k 
\end{equation}
The complex $C_{\bullet}(A)$ is the total complex of
the double complex with the differential $b + \delta$.

Let $u$ be a formal variable of degree two. The complex $(\Ccc [[u]], b+\delta + uB)$ is called {\it {the negative cyclic complex}} of $A$.

One can define explicitly a product 
\begin{equation}  \label{eq:sh}
\operatorname{sh}: \Ccc \otimes \Ccc \rightarrow \Ccc
\end{equation}
and its extension
\begin{equation}  \label{eq:sh'}
\operatorname{sh}+u\operatorname{sh}': \Ccc[[u]] \otimes \Ccc[[u]] \rightarrow \Ccc[[u]]
\end{equation}
\cite{L}. When $A$ is commutative, these are morphisms of complexes.
\subsection{Pairings between chains and cochains} \label{pairings} 
For a graded algebra $A$, for $D \in C^d(A,A)$, define
\begin{equation} \label{eq: i}
i_D(a_0 \otimes \ldots \otimes a_n)= (-1)^{|D||a_0|}a_0 D(a_1, \ldots, a_d) \otimes a_{d+1} \otimes \ldots \otimes a_n
\end{equation}
\begin{proposition}  \label{prop: properties of i}
$$ [b,i_D] = i_{\delta D};\;i_D i_E = (-1)^{|D||E|} i_{E \smile D} $$
\end{proposition} 
Now, put 
\begin{equation} \label{eq: L} 
L_D(a_0 \otimes \ldots \otimes a_n)=\sum _{k=1}^{n-d} \epsilon _k a_0 \otimes \ldots \otimes D(a_{k+1}, \ldots, a_{k+d}) \otimes \ldots \otimes a_n +
\end{equation}
$$ \sum _{k=n+1 -d}^{n} \eta _k D (a_{k+1}, \ldots, a_n, a_0, \ldots ) \otimes \ldots \otimes a_k
$$

(The second sum in the above formula is taken over all cyclic permutations such that $a_0$ is inside $D$). The signs are given by
$$ \epsilon _k = (|D| + 1)(|a_0|+\sum _{i=1}^{k} (|a_i| +1))$$
and 
$$
\eta _k = |D|+ \sum_{i \leq k}(|a_i|+1)\sum_{i \geq k}(|a_i|+1)
$$
\begin{proposition}
$$[L_D, L_E]=L_{[D,E]};\;
[b, L_D] + L_{\delta D} = 0;\;
[L_D, B] = 0$$
\end{proposition}
Now let us extend the above operations to the cyclic complex. Define
$$
S_D(a_0 \otimes \ldots \otimes a_n)= \sum_ {j\geq 0;\; k\geq j+d}\epsilon _{jk} 1  \otimes a_{k+1} \otimes \ldots a_0 \otimes \ldots 
\otimes D(a_{j+1}, \ldots, a_{j+d}) \otimes 
\ldots \otimes a_k$$
(The sum is taken over all cyclic permutations for which $a_0$ appears to the left of $D$). The signs are as follows: 
$$ \epsilon _{jk} =  (|D|+1)(\sum _{i=k+1}^{n}(|a_i|+1)+|a_0| + \sum_{i=1}^j(|a_i|+1))$$

As we will see later, all the above operations are partial cases of a unified algebraic structure for chains and cochains, cf. \ref{ss:a-infty-mod}; the sign rule for this unified construction was explained in \ref{a-infty}.

\begin{proposition} \label{prop: reinhart}
(\cite{R})
$$[b+uB, i_D + uS_D] - i_{\delta D} - uS_{\delta D} = L_D $$
\end{proposition}
The following statement implies that the differential graded Lie algebra $H^{\bullet + 1}(A) [u, \epsilon]$ with the differential $u\frac{\partial}{\partial \epsilon}$ acts on the negative cyclic homology $\HCnA$. The extension of this action to the level of cochains will me the main result of this paper.
\begin{proposition} \label{prop: gdt}
(\cite{DGT})
There exists a linear transformation $T(D,E)$ of the Hochschild chain complex, bilinear in $D, \;E \in C^{\bullet}(A,A)$, such that
\begin{eqnarray*}
[b+uB, T(D,E)] - T(\delta D, E) - (-1)^{|D|} T(D, \delta E) =\\
=[L_D , i_E + u S_E] - (-1)^{|D|+ 1}(i_{[D,E]}+uS_{[D,E]})
\end{eqnarray*}
\end{proposition}
\section{The module structure on the negative cyclic complex}\label{s:The module structure on the negative cyclic complex}
\subsection{Definitions}\label{ss:Definitions} For a monomial $Y=D_1\ldots D_n$ in $U(\gA)$, set
\begin{equation}\label{eq:Y bar}
{\overline Y}=(\ldots((D_1\circ D_2)\circ D_3)\ldots \circ D_n)\in \CochA
\end{equation}
By linearity, extend this to a map $U(\gA)\to \CochA$. It is easy to see, using induction on $n$ and Proposition \ref{prop:brace structure}, that this map is well-defined \cite{DT}.

Identify $S(\gA)$ with $U(\gA)$ as coalgebras via the Poincar\'{e}-Birkhoff-Witt map. The augmentation ideals $S(\gA)^+$ and $U(\gA)^+$ also get identified. By \begin{equation}\label{eq:coproduct}
Y\mapsto \sum Y^+_1\otimes \ldots \otimes Y_n^+
\end{equation}
denote the map
\begin{equation}\label{eq:coproduct 1}
S(\gA)^+ \to (S(\gA)^+)^{\otimes n}
\end{equation}
defined as the n-fold coproduct, followed by the $n$th power of the projection from $S(\gA)$ to $S(\gA)^+$ along $K\cdot 1$. Similarly for $U(\gA)$.
\begin{definition}\label{dfn:i and S}
For $Y\in S(\gA)^+$, define
$$
i_Y(a_0 \otimes \ldots \otimes a_n)=(-1)^{|a_0||Y|} a_0 {\overline Y}(a_1, \ldots, a_k)\otimes a_{k+1} \otimes \ldots \otimes a_n;
$$
\begin{eqnarray*}
S_Y(a_0 \otimes \ldots \otimes a_n)=\\
\sum_{n\geq 1} \sum_{i, j_1, \ldots, j_n}\pm 1 \otimes  a_{k+1}& \otimes &\ldots \otimes a_n \otimes a_0 \otimes \ldots \otimes {\overline Y}_1^+(a_{j_1}, \ldots)\otimes \ldots \otimes \overline{Y}_m^+(a_{j_m}, \ldots)\otimes \ldots \otimes a_k
\end{eqnarray*}
where the sign is 
$$ \sum_{p=1}^m  (|Y_p^+|+1)(\sum _{i=k+1}^{n}(|a_i|+1)+|a_0| + \sum_{i=1}^{j_p-1}(|a_i|+1))$$
For $Y\in S(\gA)^+$ and $D\in \gA$, define
\begin{eqnarray*}
T(D, Y)(a_0 \otimes \ldots \otimes a_n)=\\
\sum_{n\geq 1,k, j_1, \ldots, j_n}\pm D( a_{k+1},\ldots, &a_0&,\ldots, {\overline Y}_1^+(a_{j_1}, \ldots), \ldots, \overline{Y}_m^+(a_{j_m}, \ldots),\ldots, a_j)\otimes a_{j+1} \otimes \ldots \otimes a_k
\end{eqnarray*}
where the sign is 
$$(|D|+1)(\sum_{p=1}^m(|Y_p^+|+1)+ \sum_{p=1}^m  (|Y_p^+|+1)(\sum _{i=k+1}^{n}(|a_i|+1)+|a_0| + \sum_{i=1}^{j_p-1}(|a_i|+1))$$
\end{definition}
Now introduce the following differential graded algebras. Let $C(\gA[u, \epsilon])$ be the standard Chevalley-Eilenberg chain complex of the DGLA $\gA[u, \epsilon]$ over the ring of scalars $K[u]$. It carries the Chevalley-Eilenberg differential $\partial$ and the differentials $\delta$ and $\partial _{\epsilon}$ induced bu the corresponding differentials on $\gA[u, \epsilon]$. Let $C_+(\gA[u, \epsilon])$ be the augmentation co-ideal, i.e. the sum of all positive exterior powers of our DGLA. As in \eqref{eq:coproduct 1} above, the comultiplication defines maps 
$$
C_+(\gA[u, \epsilon])\mapsto C_+(\gA[u, \epsilon])^{\otimes n};
$$
$$
c\mapsto \sum c_1^+\otimes \ldots \otimes c_n^+
$$
\begin{definition}\label{dfn:Bar}
Define the associative DGLA $B(\gA[u, \epsilon])$ over $K[[u]]$ as the tensor algebra of $C_+(\gA[u, \epsilon])$ with the differential $d$ determined by 
$$dc = (\delta + \partial)c-\frac{1}{2}\sum (-1)^{|c_1^+|}c_1^+c_2^+ + u\partial _{\epsilon}c.$$ 
\end{definition}
\begin{definition}\label{dfn:Bar tw}
Let the associative DGA $B^{\operatorname {tw}}(\gA[u, \epsilon])$ over $K[[u]]$ be the tensor algebra of $C_+(\gA[u, \epsilon])$ with the differential $d$ determined by 
$$dc = (\delta + \partial)c -\frac{1}{2}\sum (-1)^{|c_1^+|}c_1^+c_2^+ +u\sum_{n=1}^{\infty} \partial _{\epsilon}c_1^+\ldots \partial _{\epsilon}c_n^+.$$ 
\end{definition}
\begin{thm} \label{thm:action of twisted bar} (cf. \cite{DT}). The following formulas define an action of the DGA $B^{\operatorname {tw}}(\gA[u, \epsilon])$ on $\CCnA$:
$$D\mapsto L_D;$$
$$\epsilon E_1 \wedge \ldots \wedge\epsilon E_n \mapsto i_Y + uS_Y$$
for $n\geq 1$;
$$ \epsilon E_1 \wedge \ldots \wedge\epsilon E_n \wedge D\mapsto T(D,Y)$$
for $n\geq 1$;
$$ \epsilon E_1 \wedge \ldots \wedge\epsilon E_n \wedge D_1 \wedge \ldots \wedge D_k\mapsto 0$$ for $k>1.$
Here $D, \; D_i, \; E_j \in \gA$ and $Y=E_1 \ldots E_n \in S(\gA)^+.$
\end{thm}
We will start the proof in section \ref{a-infty} below by recalling the $\Ai$ structure from \cite{T}, \cite{TT}. Then, in section \ref{s:proof of twisted theorem}, we will re-write the definitions in term of this $\Ai$ structure. The proof will follow from the definition of an $\Ai$ module.
\section{The $\Ai$ algebra $\CCc$}  \label{a-infty}

In this section we will construct an $\Ai$ algebra structure on the negative cyclic complex of the DGA of Hochschild cochains of any algebra $A$. The negative cyclic complex of $A$ itself will be a right $\Ai$ module over the above $\Ai$ algebra. Our construction is a direct generalization of the construction of Getzler and Jones \cite{GJ} who constructed an $\Ai$ structure on the negative cyclic complex of any commutative algebra $C$. We adapt their definition to the case when $C$ is {\em a brace algebra}, in particular the Hochschild cochain complex. 

Note that all our constructions can be carried out for a unital $\Ai$ algebra $A$. The Hochschild and cyclic complexes of $\Ai$ algebras are introduced in \cite{GJ}; as shown in \cite{G}, the Hochschild cochain complex becomes an $\Ai$ algebra; all the formulas in this section are good for the more general case. In fact they are easier to write using the $\Ai$ language, even if $A$ is a usual algebra.

Recall \cite{LS}, \cite{St} that an $A_{\infty}$ algebra is a graded vector space ${\cal{C}}$ together with a Hochschild cochain $m$ of total degree $1$, 
$$ m = \sum _{n=1}^{\infty} m_n$$
where $m_n \in C^n({\cal{C}})$ and 
$$[m,m] = 0$$

Consider the Hochschild cochain complex of a graded algebra $A$ as a differential graded associative algebra $(C^{\bullet}(A), \; \smile, \; \delta)$. Consider the Hochschild {\it chain} complex of this differential graded algebra. The total differential in this complex is $b+ \delta$; the degree of a chain is given by 
$$|D_0 \otimes \ldots \otimes D_n | = |D_0| + \sum _{i=1}^{n}(|D_i|+1)$$
where $D_i$ are Hochschild cochains.

The complex $C_{\bullet}(C^{\bullet}(A))$ contains the Hochschild cochain complex $C^{\bullet}(A)$ as a subcomplex (of zero-chains) and has the Hochschild chain complex $C_{\bullet}(A)$ as a quotient complex:
$$C^{\bullet}(A) \stackrel {i}{\longrightarrow} C_{\bullet}(C^{\bullet}(A))\stackrel {\pi}{\longrightarrow}C_{\bullet}(A)$$
(this sequence is not by any means exact). The projection on the right splits if $A$ is commutative. If not, $C_{\bullet}(A)$ is naturally a graded subspace but not a subcomplex.

\begin{thm} \label{thm: a-infty}
There is an $A_{\infty}$ structure ${\bf m}$ on $C_{\bullet}(C^{\bullet}(A))[[u]]$ such that:
\begin{itemize}
\item All ${\bf m}_n$ are $k[[u]]$-linear, $(u)$-adically continuous
\item ${\bf m}_1 = b+\delta +uB$

For $x,\;y \in C_{\bullet}(A)$:
\item $(-1)^{|x|}{\bf m}_2 (x,y) = (\operatorname{sh} + u \operatorname{sh}')(x,y)$

For $D,\;E \in C^{\bullet}(A)$:
\item $(-1)^{|D|}{\bf m}_2(D,E) = D \smile E$
\item ${\bf m}_2(1\otimes D, \; 1 \otimes E) + (-1)^{|D||E|}{\bf m}_2(1\otimes E, \; 1 \otimes D) = (-1)^{|D|}1 \otimes [D,\; E]$
\item ${\bf m}_2( D, \; 1 \otimes E) + (-1)^{(|D|+1)|E|}{\bf m}_2(1\otimes E, \;D) = (-1)^{|D|+1}[D,\; E]$
\end{itemize}
\end{thm}

Here is an explicit description of the above $A_{\infty}$ structure. We define for $n \geq 2$ 
$${\bf m}_n = {\bf m}_n^{(1)} + u{\bf m}_n^{(2)}$$
where, for 
$$a^{(k)} = D_0^{(k)} \otimes \ldots \otimes D_{N_k}^{(k)}, $$

$${\bf m}_n^{(1)}(a^{(1)}, \ldots, a^{(n)}) =  \sum \pm 
{ m}_k \{{\underline{\ldots}} , D^{(0)}_0\{{\underline{\ldots}}\} , {\underline{\ldots}} , D^{(n)}_0 \{{\underline{\ldots}}\}{\underline{\ldots}}\}\otimes {\underline{\ldots}}
$$
The space designated by $\;\underline \;$ is filled with  $D_i^{(j)},\;i>0,$ in such a way that:
\begin{itemize}
\item the cyclic order of each group $D_0 ^{(k)}, \ldots, D_{N_k}^{(k)}$ is preserved;
\item any cochain $D^{(i)}_j$ may contain some of its neighbors on the right inside the braces, provided that all of these neighbors are of the form $D^{(p)}_q$ with $p < i$.
The sign convention: any permutation contributes to the sign; the parity of $D^{(i)}_j$ is always $|D^{(i)}_j|+1 $.
\end{itemize}

$${\bf m}_n^{(2)}(a^{(1)}, \ldots, a^{(n)}) =  \sum \pm 1 \otimes 
{\underline{\ldots}} \otimes D^{(0)}_0\{{\underline{\ldots}}\} \otimes {\underline{\ldots}} \otimes D^{(n)}_0\{{\underline{\ldots}}\} \otimes {\underline{\ldots}}
$$
The space designated by $\;\underline \;$ is filled with  $D_i^{(j)},\;i>0,$ in such a way that:
\begin{itemize}
\item the cyclic order of each group $D_0 ^{(k)}, \ldots, D_{N_k}^{(k)}$ is preserved;
\item any cochain $D^{(i)}_j$ may contain some of its neighbors on the right inside the braces, provided that all of these neighbors are of the form $D^{(p)}_q$ with $p < i$.
The sign convention: any permutation contributes to the sign; the parity of $D^{(i)}_j$ is always $|D^{(i)}_j|+1 $.
\end{itemize}

\begin{remark} \label{rmk:hj} Let $A$ be a commutative algebra. Then $C_{\bullet}(A)[[u]]$ is not only a subcomplex but an $A_{\infty}$ subalgebra of $C_{\bullet}(C^{\bullet}(A))[[u]]$. The $A_{\infty}$ structure on $C_{\bullet}(A)[[u]]$ is the one from \cite{GJ}. 
\end{remark}
{\bf Proof of the Theorem.} First let us prove that ${\bf m}^{(1)}$ is an $\Ai$ structure on $C_{\bullet}(\CochA)$. Decompose it into the sum $\delta + {\widetilde {\bf m}^{(1)}}$ where $\delta$ is the differential induced by the differential on $\CochA$. We want to prove that $[\delta, {\widetilde {\bf m}^{(1)}}]+\frac{1}{2}[{\widetilde {\bf m}^{(1)}},{\widetilde {\bf m}^{(1)}}]=0$. We first compute $\frac{1}{2}[{\widetilde {\bf m}^{(1)}},{\widetilde {\bf m}^{(1)}}]$. It consists of the following terms:

(1) $ m\{\ldots D_0^{(1)} \ldots m\{ \ldots D_0^{(i+1)} \ldots D_0^{(j)}\ldots\} \ldots D_0^{(n)} \ldots\}\otimes \ldots$

where the only elements allowed inside the inner $m\{\ldots\}$ are $D_p^{(q)}$ with $i+1\leq q\leq j;$

(2) $ m\{\ldots D_0^{(1)} \ldots m\{ \ldots \} \ldots D_0^{(n)} \ldots\}\otimes \ldots$

where the only elements allowed inside the inner $m\{\ldots\}$ are $D_p^{(q)}$ for one and only $q$ (these are the contributions of the term ${\widetilde {\bf m}^{(1)}}(a^{(1)}, \ldots, ba^{(q)}, \ldots, a^{(n)}$);

(3) $ m\{\ldots D_0^{(1)} \ldots D_0^{(n)} \ldots\} \otimes \ldots \otimes  m\{ \ldots \} \otimes \ldots$

with the only requirement that the second $m\{\ldots\}$ should contain elements $D_p^{(q)}$ and $D_{p'}^{(q')}$ with $q\neq q'$.(The terms in which the second $m\{\ldots\}$ contains $D_p^{(q)}$ where all $q$'s are the same cancel out: they enter twice, as contributions from $b{\widetilde {\bf m}^{(1)}}(a^{(1)}, \ldots, a^{(q)}, \ldots, a^{(n)}$ and from ${\widetilde {\bf m}^{(1)}}(a^{(1)}, \ldots, ba^{(1)}, \ldots, a^{(n)}$). 

The collections of terms (1) and (2) differ from 

(0) $ \frac{1}{2}[m,m]\{\ldots D_0^{(1)} \ldots \ldots D_0^{(n)} \ldots\}\otimes \ldots$

by the sum of all the following terms:

$(1')$ terms as in (1), but with a requirement that in the inside $m\{\ldots \}$ an element $D_p^{(q)}$ must me present such that $q\leq i$ or $q>j;$

$(2')$ terms as in (1), but with a requirement that the inside $m\{\ldots \}$ must contain elements $D_p^{(q)}$ and $D_{p'}^{(q')}$ with $q\neq q'$.

Assume for a moment that $D_p^{(q)}$ are elements of a commutative algebra (or, more generally, of a $C_{\infty}$ algebra, i.e. a homotopy commutative algebra). Then there is no $\delta$ and ${\widetilde {\bf m}^{(1)}}={{\bf m}^{(1)}}.$ But the terms $(1')$ and $(2')$ all cancel out, as well as (3). Indeed, they all involve $m\{\ldots\}$ with some shuffles inside, and $m$ is zero on all shuffles. (the last statement is obvious for a commutative algebra, and is exactly the definition of a $C_{\infty}$ algebra).

Now, we are in a more complex situation where $D_p^{(q)}$ are Hochschild cochains (or, more generally, elements of a {\em brace algebra}). Recall that all the formulas above assume that cochains $D_p^{(q)}$ may contain their neighbors on the right inside the braces. We claim that 

(A) the terms ($1'$), $(2')$ and (3), together with (0),  cancel out with the terms constituting $[\delta, {\widetilde {\bf m}^{(1)}}]$. 

To see this, recall from \cite{KS1} the following description of brace operations. To any rooted planar tree with marked vertices one can associate an operation on Hochschild cochains. The operation
$$D\{\ldots E_1 \{\ldots \{Z_{1,1}, \ldots, Z_{1,k_1}\}, \ldots\}\ldots E_n \{\ldots \{Z_{n,1}, \ldots, Z_{n,k_n}\}\ldots\} \ldots \}$$
corresponds to a tree where $D$ is at the root, $E_i$ are connected to $D$ by edges, and so on, with $Z_{ij}$ being external vertices. The edge connecting $D$ to $E_i$ is to the left from the edge connecting $D$ to $E_j$ for $i<j$, etc. Furthermore, one is allowed to replace some of the cochains $D$, $E_i$, etc. by the cochain $m$ defining the $\Ai$ structure. In this case we leave the vertex unmarked, and regard the result as an operation whose input are cochains marking the remaining vertices (at least one vertex should remain marked). 

For a planar rooted tree $T$ with marked vertices, denote the corresponding operation by ${\bf O}_T$. The following corollary from Proposition \ref{prop:brace structure} was proven in \cite{KS1}:
$$[\delta, {\bf O}_T]=\sum_{T'} \pm {\bf O}_{T'}$$
where $T'$ are all the trees from which $T$ can be obtained by contracting an edge. One of the vertices of this new edge of $T'$ inherits the marking from the vertex to which it gets contracted; the other vertex of that edge remains unmarked. There is one restriction: the unmarked vertex of $T'$ must have more than one outgoing edge. Using this description, it is easy to see that the claim (A) is true. 

Now let us prove that
$$[\delta, {\widetilde {\bf m}^{(2)}}]+{\widetilde {\bf m}^{(1)}}\circ{{\bf m}^{(2)}}+{{\bf m}^{(2)}}\circ{\widetilde {\bf m}^{(1)}}=0$$
The summand ${{\bf m}^{(2)}}\circ{\widetilde {\bf m}^{(1)}}$ contributes both terms

(1) $D^{(1)}_0\otimes \dots \otimes D^{(2)}_0 \otimes \dots \otimes D^{(n)}_0\otimes \dots$

(2) $D^{(n)}_0\otimes \dots \otimes D^{(1)}_0 \otimes \dots \otimes D^{(n-1)}_0\otimes \dots$

twice, causing them to cancel out. Indeed, $b{\bf m}^{(2)}(a^{(1)}, \ldots, a^{(1)})$ contributes both (1) and (2); ${\widetilde{\bf m}}^{(1)}(a^{(1)},{\bf m}^{(2)}(a^{(2)}, \ldots, a^{(n)}))$ contributes (1), and ${\widetilde{\bf m}}^{(1)}({\bf m}^{(2)}(a^{(1)}, \ldots, a^{(n-1)}), a^{(n)})$ contributes (2).

(3) $1\otimes \ldots D^{(1)}_0\otimes \dots \otimes m\{ D^{(i+1)}_0 \dots D^{(j)}_0\} \otimes \ldots \otimes D^{(n)}_0\otimes \dots$

where $j\geq i$. The summand ${\widetilde {\bf m}^{(1)}}\circ{{\bf m}^{(2)}}$ consists of terms

(4) Same as (3), but with the only elements allowed inside the $m\{\ldots\}$ being $D_p^{(q)}$ with $i+1\leq q\leq j;$

(5) $1\otimes \ldots D^{(1)}_0\otimes \dots \otimes m\{ \ldots \} \otimes \ldots \otimes D^{(n)}_0\otimes \dots$

where the only elements allowed inside the $m\{\ldots\}$ are $D_p^{(q)}$ for one and only $q$. The sum of the terms (3), (4), (5) is equal to zero by the same reasoning as in the end of the proof of $[{\widetilde {\bf m}^{(1)}},{\widetilde {\bf m}^{(1)}}]=0$. 

\subsection{The $\Ai$ module structure on Hochschild chains} \label{ss:a-infty-mod}
Recall the definition of $A_{\infty}$ modules over $A_{\infty}$ algebras. First, note that for a graded space ${\cal{M}}$, the Gerstenhaber bracket $[\;,\;]$ can be extended to the space
$$\operatorname{Hom} (\overline {{\cal{C}}}^{\otimes {\bullet}}, {\cal{C}}) \oplus \operatorname{Hom} ({\cal{M}} \otimes \overline {{\cal{C}}}^{\otimes {\bullet}}, {\cal{M}})
$$

For a graded $k$-module ${\cal{M}}$, a structure of an $A_{\infty}$ module over an $A_{\infty}$ algebra ${\cal{C}}$ on ${\cal{M}}$ is a cochain 
$$ \mu = \sum _{n=1}^{\infty} \mu _n$$
$$ \mu _ n \in {Hom} ({\cal{M}} \otimes \overline {{\cal{C}}}^{\otimes {n-1}}, {\cal{M}})$$
such that 
$$[m+\mu, m+\mu ] =0$$
\begin{thm} \label{thm: a-infty mod}
On $C_{\bullet}(A)[[u]]$, there exists a structure of an $A_{\infty}$ module over the $A_{\infty}$ algebra $C_{\bullet}(C^{\bullet}(A))[[u]]$ such that: 
\begin{itemize}
\item All $\mu_n$ are $k[[u]]$-linear, $(u)$-adically continuous
\item $\mu_1 = b +uB$ on $C_{\bullet}(A)[[u]]$

For $a \in C_{\bullet}(A)[[u]]$:
\item $\mu_2 (a, D) = (-1)^{|a||D| + |a|}(i_D + u S_D)a$
\item $\mu_2 (a, 1 \otimes D) = (-1)^{|a||D|}L_D a $

For $a, \; x \in C_{\bullet}(A)[[u]]$:
$(-1)^{|a|}\mu_2(a,x) = (\operatorname{sh} + u \operatorname{sh}')(a,x)$
\end{itemize}
\end{thm}

 To obtain formulas for the structure of an $A_{\infty}$ module from Theorem
 \ref{thm: a-infty mod}, one has to assume that, in the formulas for the $\Ai$
 structure from Theorem \ref{thm: a-infty}, all $D^{(1)}_j$ are elements of
 $A$; then one has to replace braces $\{\;\}$ by the usual parentheses $(\;)$ symbolizing evaluation of a multi-linear map at elements of $A$. The proof is identical to the one for the $\Ai$ algebra case.
\section{Proof of Theorem \ref{thm:action of twisted bar}} \label{s:proof of twisted theorem} We start with two key properties of the $\Ai$ structures from section \ref{a-infty}.
\begin{proposition}\label{prop:key property 1}
Both ${\bf m}_k(c_1, \ldots, c_k)$ and $\mu _k(c_1, \ldots, c_k)$ are equal to zero if one of the arguments $c_i$, $i<k$, is of the form $1\otimes \ldots$.
\end{proposition}
\begin{proposition}\label{prop:key property 2}
For $D_i\in \gA, \, 1\leq i\geq N$, let $Y=D_1\ldots D_N \in U(\gA)$. For the $\Ai$ algebra from Theorem \ref{thm: a-infty}, put 
$$x\bullet y=(-1)^{|x|}{\bf m}_2(x,y).$$
Then
$$(1\otimes D_1)\bullet \ldots \bullet (1\otimes D_N)=\sum_{n\geq 1}1\otimes {\overline{Y}}_1^+\otimes \ldots \otimes {\overline{Y}}_n^+$$
\end{proposition}
By virtue of Proposition \ref{prop:key property 1}, the order of parentheses in the left hand side of the above formula is irrelevant.

Proposition \ref{prop:key property 1} follows immediately from the definitions, Proposition \ref{prop:key property 2} can be easily obtained by induction on $N$.

Now let us rewrite the operators from Theorem \ref{thm:action of twisted bar} in terms of the $\Ai$ structures. We replace the left module by a right module by the usual rule $x\cdot a=(-1)^{|a||x|}a\cdot x.$

For $n\geq 1$,
$$x\cdot(\epsilon E_1 \wedge \ldots \wedge\epsilon E_n)=\sum_{n\geq 1}(-1)^{|x|}\mu _{n+1}(x, \, {\overline{Y}}_1^+, \ldots , {\overline{Y}}_n^+); $$
for $n\geq 0$,
$$x\cdot(\epsilon E_1 \wedge \ldots \wedge\epsilon E_n \wedge D)= \sum_{n\geq 1}(-1)^{|x|}\mu _{n+2}(x, \, {\overline{Y}}_1^+, \ldots , {\overline{Y}}_n^+, 1\otimes D);$$
$$x\cdot(  \epsilon E_1 \wedge \ldots \wedge\epsilon E_n\wedge D_1 \wedge \ldots \wedge D_k)= 0$$ for $k>1.$
Here $D, \; D_i, \; E_j \in \gA$ and $Y=E_1 \ldots E_n \in S(\gA)^+.$

\begin{lemma} \label{lemma:Y bar as A infty morphism}
For $Y\in S(\gA)^+$ and $D\in \gA$, 
$${\overline{({\operatorname{ad}}_D Y)}}=\sum D\{{\overline{Y}}_1^+,{\overline{Y}}_2^+, \ldots, {\overline{Y}}_n^+ \}-(-1)^{(|D|+1)(|Y|+1)}{\overline Y}\{D\}.$$
In particular, for $Y\in S(\gA)^+$,
$${\overline{(\delta Y)}}=\delta {\overline{Y}}+\sum m_2\{{\overline{Y}}_1^+, {\overline{Y}}_2^+\}= \delta {\overline{Y}}+\sum (-1)^{|{Y_1^+}|+1}({\overline{Y}}_1^+\smile {\overline{Y}}_2^+)$$
\end{lemma}
Indeed, let $Y=E_1\ldots E_n.$ For $n=2$, the lemma follows from Proposition \ref{prop:brace structure}; it general, it is obtained from the same proposition by induction on $n$.
 
To prove the theorem, we have to show that 
\begin{equation}\label{eq:right module over twisted bar}
-(b+uB)(x\cdot c)+((b+uB)x)\cdot c)+(-1)^{|x|}x\cdot\partial c+(-1)^{|x|}x\cdot\delta c+
\end{equation}
\begin{equation*}
+\sum(-1)^{|x|+|c_1^+|+1}x\cdot c_1^+\cdot c_2^++(-1)^{|x|}u\sum\frac{1}{n!}x\cdot \partial_{\epsilon}c_1\cdot \ldots \cdot \partial_{\epsilon}c_n=0
\end{equation*}
Let us start by applying the $\Ai$ identity 
$$\sum \sum _{0\leq i \leq n} \pm \mu _1(\mu _{n+3}(x, {\overline{Y}}_1^+, \ldots ,{\overline{Y}}_i^+ , 1\otimes D, {\overline{Y}}_{i+1}^+ , \ldots, {\overline{Y}}_n^+, T)))+\ldots=0
$$
where $T=1\otimes F_1\otimes \ldots\otimes F_m$ is a cycle with respect to $b$ (and, automatically, to $B$). 
By virtue of Proposition \ref{prop:key property 1}, all the terms containing $1\otimes D$ in the middle vanish. The only surviving terms produce the identity
\begin{eqnarray*}\label{eq:identity 1}
\sum \pm \mu _{n-i+1}(\mu _{i+2}(x, {\overline{Y}}_1^+, \ldots ,{\overline{Y}}_i^+ , 1\otimes D), {\overline{Y}}_{i+1}^+ , \ldots, {\overline{Y}}_n^+, T))\\
+\sum \pm \mu_{n+1} (x, {\overline{Y}}_1^+, \ldots ,{\overline{Y}}_i^+\{D\} ,\ldots, {\overline{Y}}_n^+, T)+\\
+\sum \pm \mu_{n+2-j} (x, {\overline{Y}}_1^+, \ldots ,D\{{\overline{Y}}_{i+1}^+, \ldots ,{\overline{Y}}_{i+j}^+   \}, \ldots, {\overline{Y}}_n^+, T)=0
\end{eqnarray*}
When $T=1\otimes F$, $F\in \gA$, we obtain, using the first part of Lemma \ref{lemma:Y bar as A infty morphism}, the identity \eqref{eq:right module over twisted bar} for $c=\epsilon E_1\wedge \ldots \wedge \epsilon E_n\wedge  D\wedge F$. An identical computation without a $T$ at the end yields \eqref{eq:right module over twisted bar} for $c=\epsilon E_1\wedge \ldots \wedge \epsilon E_n\wedge D$. Now apply the $\Ai$ identity 
$$\sum  \pm \mu _1(\mu _{n+1}(x, {\overline{Y}}_1^+, \ldots , {\overline{Y}}_n^+))+\ldots=0
$$
We obtain
$$\sum  \pm \mu _1(\mu _{n+1}(x, {\overline{Y}}_1^+, \ldots , {\overline{Y}}_n^+))+
\sum  \pm \mu _{n+1}(\mu _1( x), {\overline{Y}}_1^+, \ldots , {\overline{Y}}_n^+)+$$
$$\sum  \pm \mu _{n+1}( x, {\overline{Y}}_1^+, \ldots ,m_1({\overline{Y}}_i^+),\ldots, {\overline{Y}}_n^+)+
\sum  \pm \mu _{n+1}( x, {\overline{Y}}_1^+, \ldots ,m_2\{{\overline{Y}}_i^+,{\overline{Y}}_{i+1}^+\},\ldots, {\overline{Y}}_n^+)+$$
$$
\sum  \pm \mu_{n-i+1} (\mu _{i+1}( x, {\overline{Y}}_1^+, \ldots ,{\overline{Y}}_i^+), {\overline{Y}}_{i+1}^+,\ldots, {\overline{Y}}_n^+)+$$
$$\sum  \pm \mu _{i+2}( x, {\overline{Y}}_1^+, \ldots ,{\overline{Y}}_i^+, 1\otimes {\overline{Y}}_{i+1}^+\otimes \ldots \otimes {\overline{Y}}_n^+)=0.$$
The first two sums in the above formula correspond to the first two terms in \eqref{eq:right module over twisted bar}; the second two sums, by virtue of the second part of Lemma \ref{lemma:Y bar as A infty morphism}, corresponds to the third term of \eqref{eq:right module over twisted bar};the fourth term of \eqref{eq:right module over twisted bar} is in our case equal to zero. The fifth sum corresponds to the fifth term of \eqref{eq:right module over twisted bar}. Now, consider the last sum in the above formula. Use Proposition \ref{prop:key property 2}, and apply the computation right after \eqref{eq:right module over twisted bar} in the case when $T=\sum 1\otimes {\overline{Y}}^+_{i+2}\otimes \ldots \otimes {\overline{Y}}^+_{n}$ and $D=1\otimes {\overline{Y}}^+_{i+1}$. Then proceed by induction on $i$. We see that the sixth sum in the formula corresponds to the sixth term of \eqref{eq:right module over twisted bar}.
\subsection{End of the proof}\label{ss:end of the proof}
It remains to pass from $B^{\operatorname {tw}}(\gA[\epsilon, u])$ to $U(\gA[\epsilon, u])$.
\begin{lemma}\label{quis of twisted and untwisted} 
The formulas 
$$D\to D;$$
$$\epsilon E_1 \wedge \ldots \wedge \epsilon E_n \mapsto \frac{1}{n!}\sum_{\sigma \in S_n}\frac{1}{n!}(\epsilon E_{\sigma_1}) E_{\sigma_2}\ldots E_{\sigma_n};$$
$$D_1\wedge \ldots D_k \wedge \epsilon E_1 \wedge \ldots \wedge \epsilon E_n\mapsto 0$$
for $k>1$ or $k=1, \, n\geq 1$
define a quasi-isomorphism of DGAs
$$B^{\operatorname {tw}}(\gA[\epsilon, u]) \to U(\gA[\epsilon, u]).$$
\end{lemma}
{\bf{Proof.}} The fact that the above map is a morphism of DGAs follows from an easy direct computation. 
To show that this is a quasi-isomorphism, consider the increasing filtration by powers of $\epsilon$. At the level of graduate quotients, $B^{\operatorname {tw}}(\gA[\epsilon, u])$ becomes the standard free resolution of $(U(\gA[\epsilon, u]), \delta)$, and the morphism is the standard map from the resolution to the algebra, therefore a quasi-isomorphism. The statement now follows from the comparison argument for spectral sequences.

To summarize, we have constructed explicitly a DGA $B^{\operatorname {tw}}(\gA[\epsilon, u])$ and the morphisms of DGAs
$$U(\gA[\epsilon, u]) \leftarrow B^{\operatorname {tw}}(\gA[\epsilon, u]) \to \operatorname{End} _{K[[u]]}(\CCnA)$$
where the morphism on the left is a quasi-isomorphism. This yields an $A_{\infty}$ morphism 
$$U(\gA[\epsilon, u])\to \operatorname{End} _{K[[u]]}(\CCnA)$$
and therefore an $L_{\infty}$ morphism 
$$\gA[\epsilon, u] \to \operatorname{End} _{K[[u]]}(\CCnA).$$

\end{document}